\documentclass[preprint,a4paper,nopreprintline]{elsarticle}
\usepackage[T1]{fontenc}
\usepackage{lmodern}
\usepackage{amsmath,amssymb,amsthm,mathrsfs}
\usepackage[margin=24mm]{geometry}
\usepackage{microtype}
\usepackage[hidelinks]{hyperref}
\biboptions{numbers,sort&compress}
\hypersetup{pdftitle={Local null controllability of a parabolic--elliptic
chemotaxis system by a Neumann boundary control},
pdfauthor={Mohmedmunavvar Mubarak Bapu}}
\allowdisplaybreaks
\journal{Nonlinear Analysis: Real World Applications}

\newtheorem{theorem}{Theorem}[section]
\newtheorem{prop}{Proposition}[section]
\newtheorem{lemma}{Lemma}[section]
\newtheorem{corollary}{Corollary}[section]
\theoremstyle{definition}
\newtheorem{definition}{Definition}[section]
\theoremstyle{remark}
\newtheorem{remark}{Remark}[section]
\numberwithin{equation}{section}

\begin{document}
\begin{frontmatter}
\title{Local null controllability of a parabolic-elliptic chemotaxis system
by a Neumann boundary control}
\author[math]{Mohmedmunavvar Mubarak Bapu\corref{cor1}}
\ead{mbapu@math.iitb.ac.in}
\cortext[cor1]{Corresponding author.}
\address[math]{Department of Mathematics, Indian Institute of Technology Bombay,
Powai, Mumbai 400076, Maharashtra, India}

\begin{abstract}
We prove local null controllability for a parabolic equation coupled to
an elliptic equation on a bounded interval. The drift depends on the
elliptic variable, and a single control acts through the Neumann condition
at the left endpoint. For every positive control time and every locally
Lipschitz drift function, sufficiently small initial data in $L^\infty$
can be driven to zero by an $L^\infty$ boundary control. To construct
the control, we extend the linear equation to a larger interval using
an even extension of the initial datum and an odd extension of the drift.
A distributed $L^2$ control in the added interval then gives a boundary
control on the original interval. A heat kernel estimate bounds the
derivative trace needed for this construction. We obtain the nonlinear
result by a set-valued fixed point argument.
\end{abstract}

\end{frontmatter}

\section{Introduction}

Let $\mathrm L>0$ and $T>0$. We consider
\begin{equation}\label{equation_one}
\begin{cases}
w_t=w_{xx}+f(u)w_x,&(t,x)\in Q:=(0,T)\times(0,\mathrm L),\\
0=u_{xx}-u-w,&(t,x)\in Q,\\
w_x(t,0)=g(t),\qquad w_x(t,\mathrm L)=0,&t\in(0,T),\\
u_x(t,0)=u_x(t,\mathrm L)=0,&t\in(0,T),\\
w(0,x)=w_0(x),&x\in(0,\mathrm L).
\end{cases}
\end{equation}
We assume throughout that $f:\mathbb R\to\mathbb R$ is locally
Lipschitz. The initial datum $w_0$ is given. Our aim is to choose the
scalar control $g$ so that $w$ vanishes at time $T$.

The elliptic equation makes the drift nonlocal in space. More precisely,
if $\partial_{xx,N}$ denotes the Neumann Laplacian, then
$$
u=-(I-\partial_{xx,N})^{-1}w.
$$
The term $f(u)w_x$ differs from the usual chemotactic flux
$\partial_x(wu_x)$. In \eqref{equation_one}, the states and controls
are real valued and have no sign restriction.

\subsection{Weak solutions and the main theorem}

We use the usual Lebesgue and Sobolev spaces on $(0,\mathrm L)$ and set
$$
V^1(Q):=
L^2(0,T;H^1(0,\mathrm L))
\cap H^1(0,T;(H^1(0,\mathrm L))^*).
$$
We equip this space with the norm
$$
\|w\|_{V^1(Q)}
:=
\|w\|_{L^2(0,T;H^1(0,\mathrm L))}
+\|w_t\|_{L^2(0,T;(H^1(0,\mathrm L))^*)}.
$$
The embedding into $C([0,T];L^2(0,\mathrm L))$ is continuous.
Endpoint values at $t=0$ and $t=T$ always refer to this continuous
representative. We write $\langle\cdot,\cdot\rangle$ for the pairing
between $(H^1)^*$ and $H^1$.

\begin{samepage}
\begin{definition}\label{def:weak-solution}
Let $w_0\in L^\infty(0,\mathrm L)$ and $g\in L^\infty(0,T)$.
A pair $(w,u)$ is a weak solution of \eqref{equation_one} if
$$
w\in V^1(Q)\cap L^\infty(Q),\qquad
u\in L^2(0,T;H^1(0,\mathrm L))\cap L^\infty(Q),
$$
$w(0)=w_0$ in $L^2(0,\mathrm L)$, and, for every
$\varphi\in L^2(0,T;H^1(0,\mathrm L))$,
\begin{align*}
&\int_0^T\langle w_t,\varphi\rangle\,dt
+\int_0^T\int_0^{\mathrm L}
\big(w_x\varphi_x-f(u)w_x\varphi\big)\,dx\,dt
+\int_0^T g(t)\varphi(t,0)\,dt=0,\\
&\int_0^T\int_0^{\mathrm L}
\big(u_x\varphi_x+u\varphi+w\varphi\big)\,dx\,dt=0.
\end{align*}
The Neumann conditions are understood through these identities.
\end{definition}
\end{samepage}

\begin{definition}
System \eqref{equation_one} is locally null controllable at time $T$
in $L^\infty(0,\mathrm L)$ if there is $\delta>0$ such that every
$w_0\in L^\infty(0,\mathrm L)$ with $\|w_0\|_{L^\infty}\leq\delta$
admits a control $g\in L^\infty(0,T)$ and a weak solution satisfying
$w(T)=0$ in $L^2(0,\mathrm L)$.
\end{definition}

Our main result is the following.
\nopagebreak

\begin{theorem}\label{main_thm}
There is a constant $c_1>0$, depending on $\mathrm L$ and $f$ but
independent of $T$, such that, for every $T>0$ and every
$w_0\in L^\infty(0,\mathrm L)$ satisfying
\begin{equation}\label{initil_cndtion}
\|w_0\|_{L^\infty(0,\mathrm L)}
\leq \exp\!\big(-c_1(1+T+1/T)\big),
\end{equation}
there is a control $g\in L^\infty(0,T)$ for which
\eqref{equation_one} has a weak solution $(w,u)$ satisfying
$w(T)=0$ in $L^2(0,\mathrm L)$.
Moreover, $u$ has a representative in
$C([0,T];H^2(0,\mathrm L))$ and $u(T)=0$.
\end{theorem}

Estimate \eqref{initil_cndtion} gives a sufficient smallness condition
for each $T>0$, with a constant independent of the control time.
We do not claim that this neighbourhood is optimal.

\subsection{Related work}

Keller and Segel \cite{18} studied models that combine diffusion with
directed movement. Horstmann \cite{horstmann2003} gives
a survey of these models and their variants. When an elliptic equation
determines an auxiliary field, the state at each time affects that
field throughout the domain. The evolution equation is then nonlocal
in space.

Fern\'andez-Cara, Limaco and de Menezes \cite{cara2013} proved
distributed null controllability for coupled parabolic and elliptic
systems using Carleman estimates and fixed point methods. For a
chemotaxis system of parabolic and elliptic type, Guo and Zhang
\cite{GuoBaozhang14} combined elliptic estimates with linear parabolic
controllability and a set-valued fixed point argument. Their later work
\cite{Guo_Zhu_Zhang} established local exact controllability to a
positive trajectory for a parabolic chemotaxis system with one
distributed control. We also use elliptic estimates and fixed point
arguments. Our boundary control construction for \eqref{equation_one}
needs one further estimate: a bound for the spatial derivative at the
interface of an extended interval.

Nonlocal couplings have also been studied in several recent works.
Wang and Zheng \cite{WangZheng2022} proved unique continuation
inequalities for a parabolic and elliptic chemotaxis system using
a priori estimates, frequency functions and localization.
Hern\'andez-Santamar\'ia and Le Balc'h \cite{HernandezLeBalch2021}
obtained local null controllability for a nonlocal semilinear heat
equation through uniform controllability of an associated reaction and
diffusion system. Bhandari and Hern\'andez-Santamar\'ia
\cite{BhandariHernandez2025} studied parabolic and elliptic systems with
additional nonlocal couplings. In the present problem, we use the
elliptic equation to obtain a bounded drift and control the parabolic
component through one endpoint.

Boundary controllability has also been studied for systems with other
couplings. Bhandari and Boyer \cite{BhandariBoyer2021} treated coupled
parabolic systems with Robin conditions and their relation to Dirichlet
controls. Bhandari, Boyer and Hern\'andez-Santamar\'ia
\cite{BhandariBoyerHernandez2021} considered systems on an interval with
Kirchhoff type conditions. For nonlinear coupling at the boundary,
Bhandari and Boyer \cite{BhandariBoyer2024} proved local exact
controllability to steady states. Bhandari, Kumbhakar and Majumdar
\cite{BhandariKumbhakarMajumdar2024} obtained local null controllability
of a nonlinear system of two parabolic equations with one Neumann
control. Boyer and
Olive \cite{BoyerOlive2024} developed a moment method for certain
multidimensional parabolic systems and showed that the controlled
boundary geometry can affect controllability. Here, the drift is
nonlocal and the control acts at only one endpoint.

For the linear step of our proof, we use the Carleman inequality of
Fern\'andez-Cara, Gonz\'alez-Burgos, Guerrero and Puel
\cite{Cara_guerrero_puel_2006} for weak solutions with nonhomogeneous
Neumann data. The same authors studied semilinear heat equations with
Fourier boundary conditions and distributed controls in
\cite{Mcaraburgospierre2006}. For general methods in parabolic
controllability and nonlinear control, we refer to
\cite{ImanuFursikovbook1996,Barbu18,coronbook07}.

\subsection{Outline of the proof}

For a prescribed bounded drift, we extend the linear equation to a
larger interval. The cited Carleman inequality gives adjoint
observability and hence, by duality, a distributed control supported
in the added interval.

To obtain a bounded boundary control, we extend the initial datum
evenly and the drift oddly. The extended solution without control is
then even. Thus only the correction with zero initial data contributes
to the derivative at the interface. Its source is supported away from
the interface. We bound this derivative using a heat kernel estimate
and a one dimensional interpolation inequality.

The linear estimate is uniform for drifts with a common $L^\infty$
bound. We then fix the source in the elliptic equation within a bounded
ball. Compactness in $L^2(Q)$ allows us to apply a set-valued fixed point
theorem. The argument applies to locally Lipschitz $f$, including
$f(0)\ne0$, and requires no endpoint compatibility conditions on the
initial datum.

Section~\ref{sc2} collects the preliminary estimates.
Section~\ref{sc3} constructs the linear boundary control using
observability and reflection. Section~\ref{sc4} proves
Theorem~\ref{main_thm} and uniqueness for a prescribed control.
Section~\ref{sc5} discusses consequences and further questions.

\section{Preliminary results}\label{sc2}

Positive constants $C$ may change between estimates and depend only on
the fixed spatial intervals. Any dependence on coefficients is stated
explicitly.
We write $\|\cdot\|_p$ for the spatial $L^p$ norm on the interval under
consideration. In particular,
$\|z(t)\|_2=\|z(t,\cdot)\|_{L^2(0,\mathrm L)}$ on $(0,\mathrm L)$.
Norms over space-time cylinders use the subscript $L^p$; for example,
$\|z\|_{L^2(Q)}^2=\int_0^T\|z(t)\|_2^2\,dt$.

\subsection{Energy spaces, endpoint values and compactness}

We recall two properties of $V^1(Q)$, valid on any bounded interval
in place of $(0,\mathrm L)$.

\begin{lemma}\label{lem:energy-compactness}
Every $z\in V^1(Q)$ has a representative in
$C([0,T];L^2(0,\mathrm L))$. For this representative,
\begin{equation}\label{energy_endpoint_bound}
\sup_{0\leq t\leq T}\|z(t)\|_2^2
\leq \frac1T\|z\|_{L^2(Q)}^2
+2\|z_t\|_{L^2(0,T;(H^1)^*)}
  \|z\|_{L^2(0,T;H^1)}.
\end{equation}
The embedding $V^1(Q)\Subset L^2(Q)$ is compact.
\end{lemma}

\begin{proof}
The Hilbert triple
$H^1(0,\mathrm L)\subset L^2(0,\mathrm L)\subset(H^1(0,\mathrm L))^*$
gives the continuous representative and the identity
\begin{equation}\label{energy_identity}
\|z(t)\|_2^2-\|z(s)\|_2^2
=2\int_s^t\langle z_t(r),z(r)\rangle\,dr,
\qquad 0\leq s\leq t\leq T;
\end{equation}
see \cite{Lions1972}. Choose $s$ where $\|z(s)\|_2^2$ is smallest;
this minimum is at most the time average. Estimating the integral by
Cauchy--Schwarz, in either time direction, gives
\eqref{energy_endpoint_bound}. Thus evaluation at any fixed time
is a bounded linear map from $V^1(Q)$ to $L^2$.

We prove compactness by a spatial approximation. Let
$(e_j)_{j\geq0}$ be an orthonormal basis of
Neumann eigenfunctions, with eigenvalues
$\lambda_j=(j\pi/\mathrm L)^2$, and let $P_N$ be the projection onto
$\operatorname{span}\{e_0,\ldots,e_N\}$. The spectral expression for
the $H^1$ norm yields
\begin{equation}\label{fourier_tail}
\|z-P_Nz\|_{L^2(Q)}^2
\leq (1+\lambda_{N+1})^{-1}
\|z\|_{L^2(0,T;H^1)}^2.
\end{equation}
If $(z_n)$ is bounded in $V^1(Q)$, each coefficient
$c_{n,j}(t)=(z_n(t),e_j)_{L^2}$ satisfies
$$
c_{n,j}'(t)=\langle (z_n)_t(t),e_j\rangle
$$
and is bounded in $H^1(0,T)$ for fixed $j$. Thus $(P_Nz_n)$ is
relatively compact in $L^2(Q)$ for each $N$. The uniform tail bound
\eqref{fourier_tail} proves total boundedness of $(z_n)$ in $L^2(Q)$.
This proves the compact embedding.
\end{proof}

\subsection{The elliptic equation}

For $\xi\in L^2(0,\mathrm L)$, consider
\begin{equation}\label{linr_ellptc}
\begin{cases}
u_{xx}-u=\xi,&x\in(0,\mathrm L),\\
u_x(0)=u_x(\mathrm L)=0.
\end{cases}
\end{equation}

\begin{prop}\label{prp21}
For $\xi\in L^2(0,\mathrm L)$, problem \eqref{linr_ellptc} has a unique
solution in $H^2(0,\mathrm L)$. More generally, for every $1<p<\infty$
and $\xi\in L^p(0,\mathrm L)$, it has a unique weak solution in
$W^{2,p}(0,\mathrm L)$, and
$$
\|u\|_{W^{2,p}(0,\mathrm L)}
\leq C_p\|\xi\|_{L^p(0,\mathrm L)}.
$$
If $\xi\in L^\infty(0,\mathrm L)$, then
\begin{equation}\label{elliptic_infty}
\|u\|_{L^\infty(0,\mathrm L)}
\leq\|\xi\|_{L^\infty(0,\mathrm L)},\qquad
\|u\|_{W^{2,\infty}(0,\mathrm L)}
\leq C\|\xi\|_{L^\infty(0,\mathrm L)}.
\end{equation}
\end{prop}

\begin{proof}
The weak equation is
$$
\int_0^{\mathrm L}(u_x\varphi_x+u\varphi)\,dx
=-\int_0^{\mathrm L}\xi\varphi\,dx,
\qquad \varphi\in H^1(0,\mathrm L).
$$
Lax--Milgram gives existence and uniqueness in $H^1$ for $\xi\in L^2$.
Since $u_{xx}=u+\xi$, we have $u\in H^2$, and the Neumann conditions
hold in the trace sense.
For the remaining estimates, we use the Green function
$$
G(x,y)=
\frac{\cosh(\min\{x,y\})\cosh(\mathrm L-\max\{x,y\})}
{\sinh\mathrm L},
\qquad
u(x)=-\int_0^{\mathrm L}G(x,y)\xi(y)\,dy.
$$
It satisfies $G\geq0$ and
$\int_0^{\mathrm L}G(x,y)\,dy=1$.
Thus $-\partial_{xx}+1$ has an inverse that is a contraction on
$L^\infty$ and, by symmetry of $G$ and interpolation, on each $L^p$.
The identity $u_{xx}=u+\xi$ and the Neumann conditions then give the
$W^{2,p}$ bound. General $L^p$ data are treated by density.
In particular, if $\xi\in L^\infty$, then
$\|u_{xx}\|_\infty\leq2\|\xi\|_\infty$ and
$\|u_x\|_\infty\leq2\mathrm L\|\xi\|_\infty$.
This proves \eqref{elliptic_infty}.
\end{proof}

Applying Proposition~\ref{prp21} to the difference of two solutions gives
\begin{equation}\label{elliptic_stability}
\|u_1-u_2\|_{H^2(0,\mathrm L)}
\leq C\|\xi_1-\xi_2\|_{L^2(0,\mathrm L)}.
\end{equation}
In one space dimension, $H^2$ embeds continuously into $C^1$, and
$H^1$ embeds continuously into $L^\infty$; see \cite{Adamas_RA_2003}.
It follows that the elliptic solution operator also maps $L^2$
continuously into $L^\infty$.

\begin{corollary}\label{cor:elliptic-time}
For $\xi\in L^2(Q)$, let $u_\xi(t)$ be the solution of
\eqref{linr_ellptc} with source $\xi(t)$. Then
\begin{equation}\label{elliptic_time_stability}
\|u_{\xi_1}-u_{\xi_2}\|_{L^2(0,T;H^2)}
\leq C\|\xi_1-\xi_2\|_{L^2(Q)}.
\end{equation}
If $\xi\in C([0,T];L^2)$, then $u_\xi\in C([0,T];H^2)$.
On the set $\{\xi: \|\xi\|_{L^\infty(Q)}\leq R\}$, the map
$\xi\mapsto f(u_\xi)$ satisfies
\begin{equation}\label{coefficient_stability}
\begin{split}
\|f(u_\xi)\|_{L^\infty(Q)}&\leq M_R:=\max_{|r|\leq R}|f(r)|,\\
\|f(u_{\xi_1})-f(u_{\xi_2})\|_{L^2(0,T;L^\infty)}
&\leq C\operatorname{Lip}_{[-R,R]}(f)
\|\xi_1-\xi_2\|_{L^2(Q)}.
\end{split}
\end{equation}
\end{corollary}

\begin{proof}
The elliptic solution operator is bounded and linear from $L^2$ to
$H^2$. It preserves strong measurability. Integrating its stability
estimate gives \eqref{elliptic_time_stability}. Its continuity also
proves the result for sources that are continuous in time.
Proposition~\ref{prp21} gives $|u_\xi|\leq R$ when $|\xi|\leq R$.
The Lipschitz estimate for $f$ on $[-R,R]$, followed by
$H^2(0,\mathrm L)\hookrightarrow L^\infty(0,\mathrm L)$, gives
\eqref{coefficient_stability}. For every finite $p>1$, bounded sources
also give $u_\xi\in L^\infty(0,T;W^{2,p})$ by the same argument.
\end{proof}

\subsection{The linear boundary equation}

Let $B\in L^\infty(Q)$. We consider
\begin{equation}\label{lnr_prblic}
\begin{cases}
w_t=w_{xx}+B(t,x)w_x,&(t,x)\in Q,\\
w_x(t,0)=g(t),\qquad w_x(t,\mathrm L)=0,&t\in(0,T),\\
w(0,x)=w_0(x),&x\in(0,\mathrm L).
\end{cases}
\end{equation}

\begin{prop}\label{prop21}
For $w_0\in L^2(0,\mathrm L)$ and $g\in L^2(0,T)$,
problem \eqref{lnr_prblic} has a unique weak solution $w\in V^1(Q)$.
With $M=\|B\|_{L^\infty(Q)}$, the solution satisfies
\begin{equation}\label{sln_estmt}
\|w\|_{V^1(Q)}
+\|w\|_{C([0,T];L^2(0,\mathrm L))}
\leq \exp\!\big(C(1+M^2)(1+T)\big)
\big(\|w_0\|_{L^2(0,\mathrm L)}+\|g\|_{L^2(0,T)}\big).
\end{equation}
If $g=0$ and $w_0\in L^\infty(0,\mathrm L)$, then
\begin{equation}\label{homogeneous_maximum}
\|w\|_{L^\infty(Q)}
\leq\|w_0\|_{L^\infty(0,\mathrm L)}.
\end{equation}
\end{prop}

\begin{proof}
For a.e.\ $t$, define
$$
a(t;h,\varphi)=
\int_0^{\mathrm L}h_x\varphi_x\,dx
-\int_0^{\mathrm L}Bh_x\varphi\,dx,
\qquad h,\varphi\in H^1(0,\mathrm L).
$$
The boundary functional is
$\ell(t)(\varphi)=-g(t)\varphi(0)$.
The trace inequality
$$
|h(0)|^2\leq C\big(\|h\|_2^2+\|h\|_2\|h_x\|_2\big)
\leq C\|h\|_{H^1}^2
$$
shows that $\ell\in L^2(0,T;(H^1)^*)$.
The form is bounded and satisfies the G{\aa}rding inequality
$$
a(t;h,h)\geq
\frac12\|h_x\|_2^2-\frac{M^2}{2}\|h\|_2^2.
$$
To obtain a coercive form, choose $\mu=(1+M^2)/2$ and set
$w(t)=\exp\!\big(\mu t\big)z(t)$. The equation for $z$ has the form
$a(t;z,\varphi)+\mu(z,\varphi)_{L^2}$, initial datum $w_0$ and
forcing $\exp\!\big(-\mu t\big)\ell(t)$. The shifted form is bounded and
coercive on $H^1$. The variational existence theorem
\cite{Lions1972} applies to $z$, and multiplication by $\exp\!\big(\mu t\big)$
gives a solution of the original equation in $V^1(Q)$.

Testing with $w$ and applying the trace bound gives
\begin{equation}\label{bfr_grnwl_2}
\frac{d}{dt}\|w(t)\|_2^2+\|w_x(t)\|_2^2
\leq C(1+M^2)\|w(t)\|_2^2+C|g(t)|^2.
\end{equation}
Gronwall's inequality bounds the $C([0,T];L^2)$ norm and the
$L^2(0,T;H^1)$ norm.
The weak equation also gives
$$
\|w_t(t)\|_{(H^1)^*}
\leq C(1+M)\|w_x(t)\|_2+C|g(t)|.
$$
These inequalities imply \eqref{sln_estmt}.
Applying them to the difference of two solutions proves uniqueness.

If $g=0$, set $k=\|w_0\|_\infty$ and test with $(w-k)_+$.
Since $w_x=\partial_x(w-k)_+$ where $w>k$, one obtains
$$
\frac{d}{dt}\|(w-k)_+(t)\|_2^2
\leq M^2\|(w-k)_+(t)\|_2^2.
$$
Its initial value is zero. Thus $w\leq k$.
The same argument for $-w$ proves \eqref{homogeneous_maximum}.
\end{proof}

\subsection{A bounded derivative trace at the interface}

Let $\ell>0$, let $I=(-\ell,\mathrm L)$, and choose an open interval
$\widetilde\omega$ whose closure is contained in $(-\ell,0)$.
Set $Q_\ell=(0,T)\times I$ and
$$
H_N^2(I):=\{h\in H^2(I):h_x(-\ell)=h_x(\mathrm L)=0\}.
$$
The next lemma gives the regularity estimate needed to construct the
boundary control. We estimate the spatial derivative at the interface,
where the control will act.

\begin{lemma}\label{lem:boundary-trace}
Let $\widetilde B\in L^\infty(Q_\ell)$ and
$v\in L^2((0,T)\times\widetilde\omega)$.
The solution of
\begin{equation}\label{trace_auxiliary}
\begin{cases}
Y_t=Y_{xx}+\widetilde B Y_x+\chi_{\widetilde\omega}v,
&(t,x)\in Q_\ell,\\
Y_x(t,-\ell)=Y_x(t,\mathrm L)=0,&t\in(0,T),\\
Y(0,x)=0,&x\in I
\end{cases}
\end{equation}
satisfies
$$
Y\in H^1(0,T;L^2(I))\cap L^2(0,T;H_N^2(I))
\cap C([0,T];H^1(I)).
$$
If $M=\|\widetilde B\|_{L^\infty(Q_\ell)}$, then
\begin{equation}\label{H2-estimate}
\begin{split}
&\|Y\|_{C([0,T];H^1(I))}
+\|Y\|_{L^2(0,T;H^2(I))}
+\|Y_t\|_{L^2(Q_\ell)}
+\|Y_x\|_{L^6(Q_\ell)}\\
&\hspace{15mm}\leq
\exp\!\big(C(1+M^2)(1+T)\big)
\|v\|_{L^2((0,T)\times\widetilde\omega)}.
\end{split}
\end{equation}
Furthermore, the interior trace $Y_x(\cdot,0)$ belongs to
$L^\infty(0,T)$ and
\begin{equation}\label{bounded_interface_trace}
\|Y_x(\cdot,0)\|_{L^\infty(0,T)}
\leq \exp\!\big(C(1+M^2)(1+T)\big)
\|v\|_{L^2((0,T)\times\widetilde\omega)}.
\end{equation}
The constants depend only on $I$ and $\widetilde\omega$.
\end{lemma}

\begin{proof}
\emph{Energy estimates.}
Testing first with $Y$ and then with $-Y_{xx}$ gives

\begin{align*}
\frac{d}{dt}\|Y\|_2^2+\|Y_x\|_2^2
&\leq C(1+M^2)\|Y\|_2^2+C\|v(t)\|_2^2,\\
\frac{d}{dt}\|Y_x\|_2^2+\|Y_{xx}\|_2^2
&\leq CM^2\|Y_x\|_2^2+C\|v(t)\|_2^2.
\end{align*}
To justify the tests, we first work with Galerkin approximations in
the Neumann eigenfunctions. Let $P_n$ project onto the first $n+1$
eigenfunctions and solve the finite system
$$
\partial_tY_n=\partial_{xx}Y_n+
P_n(\widetilde B\partial_xY_n+\chi_{\widetilde\omega}v),
\qquad Y_n(0)=0.
$$
Both $Y_n$ and $-\partial_{xx}Y_n$ are admissible tests in this
system. Orthogonality removes the projection from the tested terms.
The estimates are therefore uniform in $n$ and require only
measurability and boundedness of $\widetilde B$. Since $Y(0)=0$,
integration and Gronwall's inequality give
\begin{align*}
\sup_{0\leq t\leq T}\|Y(t)\|_2^2
+\int_0^T\|Y_x(t)\|_2^2\,dt
&\leq Ce^{C(1+M^2)T}\|v\|_{L^2}^2,\\
\sup_{0\leq t\leq T}\|Y_x(t)\|_2^2
+\int_0^T\|Y_{xx}(t)\|_2^2\,dt
&\leq Ce^{CM^2T}\|v\|_{L^2}^2.
\end{align*}
By weak compactness, a subsequence converges to a limit. Uniqueness
shows that this limit is the energy solution. The estimates bound both the
$L^2(0,T;H^2)$ norm and the $L^\infty(0,T;H^1)$ norm.
The equation then gives
$$
\|Y_t\|_{L^2(Q_\ell)}
\leq\|Y_{xx}\|_{L^2(Q_\ell)}
+M\|Y_x\|_{L^2(Q_\ell)}+\|v\|_{L^2}.
$$
The continuity into $H^1$ follows from
$Y\in H^1(0,T;L^2)\cap L^2(0,T;H_N^2)$
\cite{Lions1972}.
No derivative of $\widetilde B$ is used.

For a.e.\ $t$, $Y_x(t)\in H_0^1(I)$.
The one dimensional Gagliardo--Nirenberg inequality gives
$$
\|Y_x(t)\|_6^6
\leq C\|Y_x(t)\|_2^4\|Y_{xx}(t)\|_2^2.
$$
Integration in time proves the $L^6$ part of \eqref{H2-estimate}.

\smallskip\noindent
\emph{Heat kernel bounds.}
Let $S(t)$ be the Neumann heat semigroup on $I$, with kernel
$K_N(t,x,y)$. We use the estimates
\begin{equation}\label{heat_gradient_six}
\|\partial_xS(r)h\|_{L^\infty(I)}
\leq C(1+r^{-7/12})\|h\|_{L^6(I)},\qquad r>0,
\end{equation}
and, with $d=\operatorname{dist}(0,\widetilde\omega)>0$,
\begin{equation}\label{heat_off_diagonal}
|\partial_xK_N(r,0,y)|
\leq
\begin{cases}
Cr^{-1}\exp(-d^2/(Cr)),&0<r\leq1,\\
C,&r\geq1,
\end{cases}
\qquad y\in\widetilde\omega.
\end{equation}
These bounds follow from the reflection formula. Put $a=-\ell$,
$D=\mathrm L+\ell$ and $G(r,z)=(4\pi r)^{-1/2}\exp\!\big(-z^2/(4r)\big)$.
Then
$$
K_N(r,x,y)=\sum_{k\in\mathbb Z}
\big[G(r,x-y+2kD)+G(r,x+y-2a+2kD)\big].
$$
For a single Gaussian,
$$
\partial_zG(r,z)=-\frac{z}{2r}G(r,z),\qquad
|\partial_zG(r,z)|\leq Cr^{-1}\exp\!\big(-z^2/(8r)\big).
$$
At $x=0$ and $y\in\widetilde\omega$, all the image distances in the
reflection formula are at least $d$. The distances of the remaining
images grow linearly with $|k|$. Summing their Gaussian bounds for
$0<r\leq1$ proves \eqref{heat_off_diagonal}.
For the estimate on the full interval, a change of variable gives
$$
\|\partial_zG(r,\cdot)\|_{L^{6/5}(\mathbb R)}
=Cr^{-1+1/(2(6/5))}=Cr^{-7/12}.
$$
The reflection formula gives the same bound for
$\partial_xK_N(r,x,\cdot)$, uniformly in $x\in I$, for $r\leq1$.
H\"older's inequality in the kernel integral proves
\eqref{heat_gradient_six} in this range. For $r\geq1$, write
$S(r)=S(1)S(r-1)$ and use the $L^6$ contraction of $S(r-1)$.
The smooth kernel at time $1$ gives the remaining uniform bounds.

\smallskip\noindent
\emph{The trace at the interface.}
The variation of constants formula for \eqref{trace_auxiliary} is
$$
Y(t)=\int_0^tS(t-s)
\big(\widetilde B(s)Y_x(s)+\chi_{\widetilde\omega}v(s)\big)\,ds.
$$
For the drift term, \eqref{heat_gradient_six} and H\"older's
inequality bound the derivative at $x=0$ by
\begin{equation}\label{trace_drift_term}
\begin{split}
&CM\int_0^t(1+(t-s)^{-7/12})\|Y_x(s)\|_6\,ds\\
&\hspace{10mm}\leq
CM(T^{5/6}+T^{1/4})\|Y_x\|_{L^6(Q_\ell)}.
\end{split}
\end{equation}
Indeed, $(7/12)(6/5)=7/10<1$.
For the second term, \eqref{heat_off_diagonal} gives
$$
\int_0^T\int_{\widetilde\omega}
|\partial_xK_N(r,0,y)|^2\,dy\,dr
\leq C\int_0^1r^{-2}\exp\!\big(-d^2/(Cr)\big)\,dr+CT
\leq C(1+T).
$$
Thus the source contribution to the derivative at the interface is
bounded uniformly in $t$ by $C(1+T)^{1/2}\|v\|_{L^2}$.
We now check that differentiating the formula gives the derivative
trace. Approximate
$\widetilde B Y_x$ in $L^6(Q_\ell)$ and $v$ in
$L^2((0,T)\times\widetilde\omega)$ by smooth functions, extending the
second approximation by zero. Let $Y_n$ be the corresponding heat
solutions with zero initial data. The $L^2$ heat regularity estimate
gives $Y_n\to Y$ in $H^1(0,T;L^2)\cap L^2(0,T;H_N^2)$.
Hence $(Y_n)_x(\cdot,0)\to Y_x(\cdot,0)$ in $L^2(0,T)$.
The kernel estimates give convergence of these traces in
$L^\infty(0,T)$ to the differentiated formula. Thus the two limits
agree almost everywhere. We combine these bounds with
\eqref{H2-estimate} and absorb the polynomial factors into the
exponential to obtain \eqref{bounded_interface_trace}.
\end{proof}

\section{Null controllability for the linear boundary equation}\label{sc3}

In this section we fix
$$
\ell=\mathrm L,\qquad
I=(-\ell,\mathrm L),\qquad
\widetilde\omega=(-3\mathrm L/4,-\mathrm L/2).
$$
For $B\in L^\infty(Q)$, define
\begin{equation}\label{vlue_kappa}
\kappa=(1+\|B\|_{L^\infty(Q)}^2)(1+T+1/T).
\end{equation}

\begin{theorem}\label{thm31}
Let $B\in L^\infty(Q)$ and $w_0\in L^2(0,\mathrm L)$.
There is a control $g\in L^\infty(0,T)$ such that the weak solution
of \eqref{lnr_prblic} satisfies $w(T)=0$ in $L^2(0,\mathrm L)$ and
\begin{equation}\label{cntrl_estmt}
\|g\|_{L^\infty(0,T)}
\leq \exp\!\big(C\kappa\big)\|w_0\|_{L^2(0,\mathrm L)}.
\end{equation}
If, in addition, $w_0\in L^\infty(0,\mathrm L)$, the control can be
chosen so that
\begin{equation}\label{linear_controlled_state}
\|w\|_{V^1(Q)}+\|w\|_{L^\infty(Q)}
+\|g\|_{L^\infty(0,T)}
\leq \exp\!\big(C\kappa\big)\|w_0\|_{L^\infty(0,\mathrm L)}.
\end{equation}
The constant $C$ depends only on $\mathrm L$.
\end{theorem}

\subsection{The adjoint equation and observability}

On the enlarged interval, consider
\begin{equation}\label{extended_system}
\begin{cases}
Z_t=Z_{xx}+\widetilde B Z_x+\chi_{\widetilde\omega}v,
&(t,x)\in Q_\ell,\\
Z_x(t,-\ell)=Z_x(t,\mathrm L)=0,&t\in(0,T),\\
Z(0,x)=Z_0(x),&x\in I.
\end{cases}
\end{equation}

The adjoint associated with \eqref{extended_system} is
\begin{equation}\label{ad_main}
\begin{cases}
-\varphi_t=\varphi_{xx}-(\widetilde B\varphi)_x,&(t,x)\in Q_\ell,\\
\varphi_x-\widetilde B\varphi=0,&x=-\ell,\mathrm L,\quad t\in(0,T),\\
\varphi(T,x)=\varphi_T(x),&x\in I.
\end{cases}
\end{equation}
For every $\varphi_T\in L^2(I)$, there is a unique solution in
$V^1(Q_\ell)$. Its weak equation is
\begin{equation}\label{adjoint_weak_form}
-\langle\varphi_t,\eta\rangle
+\int_I\varphi_x\eta_x\,dx
-\int_I\widetilde B\varphi\eta_x\,dx=0,
\qquad \eta\in H^1(I),
\end{equation}
and $\varphi(T)=\varphi_T$. Time reversal and the variational argument
of Proposition~\ref{prop21}, applied to the transposed form, give
existence and uniqueness. The associated energy estimate is
\begin{equation}\label{adjoint_energy}
-\frac{d}{dt}\|\varphi(t)\|_2^2+\|\varphi_x(t)\|_2^2
\leq M^2\|\varphi(t)\|_2^2,
\qquad M=\|\widetilde B\|_{L^\infty(Q_\ell)}.
\end{equation}

\begin{remark}\label{rem:adjoint}
The boundary data in \eqref{ad_main} prescribe the flux
$\varphi_x-\widetilde B\varphi$. A bounded measurable coefficient
$\widetilde B$ need not have endpoint traces. The weak equation
\eqref{adjoint_weak_form} states the flux condition using only
values of the coefficient inside the interval. We use this formulation
in the Carleman estimate below.
\end{remark}

We recall the Carleman estimate used in the proof. Choose an
open interval $\omega'$ with
$\overline{\omega'}\subset\widetilde\omega$ and a function
$\alpha\in C^2(\overline I)$ such that
$$
\alpha>0\text{ in }I,\qquad
\alpha(-\ell)=\alpha(\mathrm L)=0,\qquad
|\alpha_x|>0\text{ on }\overline I\setminus\omega'.
$$
Such a function can be chosen to have its unique critical point in
$\omega'$. Define the positive weights
\begin{equation}\label{carleman_weights}
\xi(t,x)=\frac{\exp\!\big(\lambda\alpha(x)\big)}{t(T-t)},
\qquad
\beta(t,x)=
\frac{\exp\!\big(2\lambda\|\alpha\|_\infty\big)-\exp\!\big(\lambda\alpha(x)\big)}{t(T-t)}.
\end{equation}
For $s,\lambda>0$, write
$$
\mathcal I_{s,\lambda}(\varphi)
=\int_{Q_\ell}\exp\!\big(-2s\beta\big)
\big(s\lambda^2\xi|\varphi_x|^2
+s^3\lambda^4\xi^3|\varphi|^2\big)\,dx\,dt.
$$

Let $\Sigma=(0,T)\times\{-\ell,\mathrm L\}$, with integration over
$\Sigma$ defined as the sum of the two endpoint integrals. The following
statement is \cite[Theorem~1]{Cara_guerrero_puel_2006} on an interval.
The weight denoted by $\alpha$ in that source is denoted here by
$\beta$, as in \eqref{carleman_weights}.

\begin{theorem}[Carleman estimate]\label{thm:adjoint-full}
Let $f_1,f_2\in L^2(Q_\ell)$, $f_3,f_4\in L^2(0,T)$ and
$\varphi_T\in L^2(I)$. Consider the weak solution of
\begin{equation}\label{eq:adjoint-full}
\begin{cases}
-\varphi_t=\varphi_{xx}+f_1+\partial_xf_2,&(t,x)\in Q_\ell,\\
(\varphi_x+f_2)(t,-\ell)=f_3(t),&t\in(0,T),\\
(\varphi_x+f_2)(t,\mathrm L)=f_4(t),&t\in(0,T),\\
\varphi(T,x)=\varphi_T(x),&x\in I.
\end{cases}
\end{equation}
There are positive constants $\bar\lambda,\sigma_1,\sigma_2,C$,
depending only on the fixed intervals and the spatial weight, such
that, for $\lambda\geq\bar\lambda$ and
$s\geq\sigma_1(\exp\!\big(\sigma_2\lambda\big)T+T^2)$,
\begin{equation}\label{eq:carleman-observability}
\begin{split}
&\mathcal I_{s,\lambda}(\varphi)
+s^2\lambda^3\int_\Sigma \exp\!\big(-2s\beta\big)\xi^2|\varphi|^2\,d\sigma\,dt\\
&\quad\leq C\Bigg[
\int_{Q_\ell}\exp\!\big(-2s\beta\big)
\big(|f_1|^2+s^2\lambda^2\xi^2|f_2|^2\big)\,dx\,dt\\
&\qquad+s\lambda\int_0^T\Big(
(\exp\!\big(-2s\beta\big)\xi)(t,-\ell)|f_3(t)|^2
+(\exp\!\big(-2s\beta\big)\xi)(t,\mathrm L)|f_4(t)|^2\Big)\,dt\\
&\qquad+s^3\lambda^4\int_0^T\int_{\widetilde\omega}
\exp\!\big(-2s\beta\big)\xi^3|\varphi|^2\,dx\,dt\Bigg].
\end{split}
\end{equation}
The boundary conditions in \eqref{eq:adjoint-full} are variational
conditions on the flux; separate traces of $f_2$ are not required.
\end{theorem}

At the left endpoint, the outward normal is $-1$, so the outward flux
data in the cited theorem are $-f_3$ and $f_4$. Their squares give the
two terms in \eqref{eq:carleman-observability}. Both endpoints are
therefore included, even when the separate terms in the flux have no
traces. The variational equation for \eqref{eq:adjoint-full} is
$$
-\langle\varphi_t,\eta\rangle+\int_I\varphi_x\eta_x\,dx
=\int_I(f_1\eta-f_2\eta_x)\,dx
+f_4(t)\eta(\mathrm L)-f_3(t)\eta(-\ell),
\qquad\eta\in H^1(I).
$$

\begin{prop}[Observability inequality]\label{prop:ob_ineq}
Every weak solution of \eqref{ad_main} satisfies
\begin{equation}\label{ob_ineq}
\|\varphi(0)\|_{L^2(I)}^2
\leq \exp\!\big(C(1+1/T+M^2+TM^2)\big)
\int_0^T\int_{\widetilde\omega}|\varphi|^2\,dx\,dt.
\end{equation}
The constant $C$ depends only on the fixed intervals.
\end{prop}

\begin{proof}
\emph{Applying the Carleman estimate.}
In Theorem~\ref{thm:adjoint-full}, take
$$
f_1=0,\qquad f_2=-\widetilde B\varphi,\qquad f_3=f_4=0.
$$
This gives exactly \eqref{ad_main}. In particular,
$f_2\in L^2(Q_\ell)$.
After dropping the nonnegative boundary term in that estimate, we obtain
\begin{equation}\label{carlmn_estmt}
\begin{split}
\mathcal I_{s,\lambda}(\varphi)
\leq {}&Cs^2\lambda^2
\int_{Q_\ell}\exp\!\big(-2s\beta\big)\xi^2
|\widetilde B\varphi|^2\,dx\,dt\\
&+Cs^3\lambda^4\int_0^T\int_{\widetilde\omega}
\exp\!\big(-2s\beta\big)\xi^3|\varphi|^2\,dx\,dt.
\end{split}
\end{equation}
Here $\lambda\geq\bar\lambda$ and
$s\geq\sigma_1(\exp\!\big(\sigma_2\lambda\big)T+T^2)$, with constants depending
only on the intervals and the chosen weight.

Fix $\lambda=\lambda_0:=\max\{1,\bar\lambda\}$ throughout the proof.
Since $\alpha\geq0$,
\begin{equation}\label{lowbd1}
\xi(t,x)\geq\frac4{T^2},\qquad (t,x)\in Q_\ell.
\end{equation}
The first term on the right of \eqref{carlmn_estmt} is bounded by
$$
\frac{CM^2T^2}{4s\lambda_0^2}
\int_{Q_\ell}\exp\!\big(-2s\beta\big)s^3\lambda_0^4\xi^3|\varphi|^2\,dx\,dt.
$$
This term can be absorbed into the left side if we take
\begin{equation}\label{carleman_parameter_choice}
s=C_0\big(T+T^2(1+M^2)\big),
\end{equation}
where the geometric constant $C_0$ is sufficiently large. This choice
also meets the lower threshold on $s$. It yields
\begin{equation}\label{using_main_carleman}
\mathcal I_{s,\lambda_0}(\varphi)
\leq Cs^3\lambda_0^4
\int_0^T\int_{\widetilde\omega}
\exp\!\big(-2s\beta\big)\xi^3|\varphi|^2\,dx\,dt.
\end{equation}
Neither $\lambda_0$ nor $C_0$ depends on $M$ or $T$.

\smallskip\noindent
\emph{Removing the weights.}
To estimate the solution at the initial time, let
$J_T=[T/4,3T/4]$. On $J_T\times I$, the fixed spatial weight
gives $\beta\leq C/T^2$ and $\xi\geq c/T^2$. Consequently,
\begin{equation}\label{central_weight_bound}
\int_{J_T}\|\varphi(t)\|_2^2\,dt
\leq C\frac{T^6}{s^3}\exp\!\big(Cs/T^2\big)
\mathcal I_{s,\lambda_0}(\varphi).
\end{equation}
On the whole cylinder,
$$
\frac{\beta}{\xi}
=\exp\!\big(2\lambda_0\|\alpha\|_\infty-\lambda_0\alpha(x)\big)-1
\geq \exp\!\big(\lambda_0\|\alpha\|_\infty\big)-1>0.
$$
The function $r\mapsto r^3\exp\!\big(-cr\big)$ is bounded for $r\geq0$.
Thus $s^3\lambda_0^4\xi^3\exp\!\big(-2s\beta\big)\leq C$, and
\eqref{using_main_carleman} implies
\begin{equation}\label{unweighted_observation}
\mathcal I_{s,\lambda_0}(\varphi)
\leq C\int_0^T\int_{\widetilde\omega}|\varphi|^2\,dx\,dt.
\end{equation}

The energy inequality \eqref{adjoint_energy} gives
$\|\varphi(0)\|_2^2\leq \exp\!\big(M^2t\big)\|\varphi(t)\|_2^2$ for
$0\leq t\leq T$. Integrating on $J_T$ yields
\begin{equation}\label{t/3to3t/4}
\|\varphi(0)\|_2^2
\leq\frac2T \exp\!\big(M^2T\big)\int_{J_T}\|\varphi(t)\|_2^2\,dt.
\end{equation}
Combining \eqref{central_weight_bound}--\eqref{t/3to3t/4} gives
$$
\|\varphi(0)\|_2^2
\leq C\frac{T^5}{s^3}\exp\!\big(Cs/T^2+M^2T\big)
\int_0^T\int_{\widetilde\omega}|\varphi|^2\,dx\,dt.
$$
By \eqref{carleman_parameter_choice},
$$
\frac{T^5}{s^3}\leq C\frac{T^2}{(1+T)^3}\leq C,
\qquad
\frac{s}{T^2}=C_0(1/T+1+M^2).
$$
This proves \eqref{ob_ineq}. All the estimates apply directly to the
weak adjoint solution, since the cited Carleman inequality allows
$L^2$ divergence sources and variational flux data.
\end{proof}

\subsection{A distributed control by duality}

\begin{prop}\label{prop:internal-control}
Let $\widetilde B\in L^\infty(Q_\ell)$ and $Z_0\in L^2(I)$.
There exists $v\in L^2((0,T)\times\widetilde\omega)$ such that the
solution of \eqref{extended_system} satisfies $Z(T)=0$ and
\begin{equation}\label{internal_control_cost}
\|v\|_{L^2((0,T)\times\widetilde\omega)}
\leq
\exp\!\big(C(1+1/T+M^2+TM^2)\big)\|Z_0\|_{L^2(I)}.
\end{equation}
\end{prop}

\begin{proof}
Let $\mathcal R$ be the linear subspace of
$L^2((0,T)\times\widetilde\omega)$ consisting of restrictions of
solutions of \eqref{ad_main}, as $\varphi_T$ ranges over $L^2(I)$.
Define a linear functional on this subspace by
$$
\mathcal F(\varphi|_{(0,T)\times\widetilde\omega})
=-(Z_0,\varphi(0))_{L^2(I)}.
$$
If two adjoint solutions agree on the observation set, observability
applied to their difference shows that their initial values agree.
Thus $\mathcal F$ is well defined. Proposition~\ref{prop:ob_ineq} also
gives
$$
|\mathcal F(\varphi|_{(0,T)\times\widetilde\omega})|
\leq \exp\!\big(C(1+1/T+M^2+TM^2)\big)\|Z_0\|_2
\|\varphi\|_{L^2((0,T)\times\widetilde\omega)}.
$$
The Hahn--Banach theorem extends $\mathcal F$ to the full observation
space with the same norm. By the Riesz representation theorem, there
is $v$ in that space satisfying \eqref{internal_control_cost} and
\begin{equation}\label{duality_control_identity}
\int_0^T\int_{\widetilde\omega}v\varphi\,dx\,dt
=-(Z_0,\varphi(0))_{L^2(I)}
\qquad\text{for every }\varphi_T\in L^2(I).
\end{equation}

Let $Z\in V^1(Q_\ell)$ solve \eqref{extended_system} with this $v$.
The weak forward and adjoint equations give
\begin{equation}\label{forward_adjoint_identity}
(Z(T),\varphi_T)_{L^2(I)}
-(Z_0,\varphi(0))_{L^2(I)}
=\int_0^T\int_{\widetilde\omega}v\varphi\,dx\,dt.
\end{equation}
To see this, apply the product rule in the energy triple to $Z$ and
$\varphi$. The weak equations cancel the diffusion and drift terms,
leaving the control integral on the right.
Equations \eqref{duality_control_identity} and
\eqref{forward_adjoint_identity} imply
$(Z(T),\varphi_T)_{L^2(I)}=0$ for every $\varphi_T\in L^2(I)$.
Hence $Z(T)=0$. This constructs the distributed null control by the
duality argument in \cite[Theorem 2]{Cara_guerrero_puel_2006}, applied
to the extended equation.
\end{proof}

\subsection{Reflection and restriction}

\begin{proof}[Proof of Theorem~\ref{thm31}]
Define the odd extension of the drift and the even extension of the
initial datum by
\begin{equation}\label{reflected_data}
\widetilde B(t,x)=
\begin{cases}
B(t,x),&0<x<\mathrm L,\\
-B(t,-x),&-\mathrm L<x<0,
\end{cases}
\qquad
Z_0(x)=w_0(|x|).
\end{equation}
The value of $\widetilde B$ at $x=0$ can be assigned arbitrarily.
In particular,
$$
\|\widetilde B\|_{L^\infty(Q_\ell)}=\|B\|_{L^\infty(Q)},
\qquad
\|Z_0\|_{L^2(I)}=\sqrt2\|w_0\|_{L^2(0,\mathrm L)}.
$$
The extension may have a jump at $0$, which is allowed in
Proposition~\ref{prop:internal-control}.

Let $v$ be the control given by that proposition and $Z$ its state.
Let $Z^0$ solve the extended equation with $v=0$ and initial datum
$Z_0$.
Since $\widetilde B$ is odd and $Z_0$ is even,
$Z^0(t,-x)$ solves the same problem as $Z^0(t,x)$.
Uniqueness gives
$$
Z^0(t,-x)=Z^0(t,x).
$$
For $\eta\in H^1(0,\mathrm L)$, its even extension belongs to $H^1(I)$.
Testing the weak equation for $Z^0$ with this extension gives even
integrands in the diffusion term and in the term containing
$\widetilde B Z_x^0$. Dividing the identity by two shows that
$w^0:=Z^0|_{(0,\mathrm L)}$ is the solution of
\eqref{lnr_prblic} with $g=0$.
This argument uses only the weak equation, so no pointwise derivative
of the initial datum is needed.

Set $Y=Z-Z^0$. It solves \eqref{trace_auxiliary} with zero initial data.
By Lemma~\ref{lem:boundary-trace}, the function
\begin{equation}\label{boundary_control_definition}
g(t):=Y_x(t,0)
\end{equation}
belongs to $L^\infty(0,T)$.
Since $\widetilde\omega\subset(-\ell,0)$, restriction to the original
interval gives
$$
w:=w^0+Y|_{(0,\mathrm L)}=Z|_{(0,\mathrm L)}.
$$
The source vanishes on the original interval, so the correction satisfies
$Y_t=Y_{xx}+BY_x$ on $(0,\mathrm L)$.
Its $H^2$ regularity and its zero derivative at $\mathrm L$ give
$$
\langle Y_t,\eta\rangle
+\int_0^{\mathrm L}Y_x\eta_x\,dx
-\int_0^{\mathrm L}BY_x\eta\,dx
=-g(t)\eta(0),\qquad \eta\in H^1(0,\mathrm L).
$$
Adding the weak equation for $w^0$ proves that $w$ solves
\eqref{lnr_prblic} with control \eqref{boundary_control_definition}.
The condition at $\mathrm L$ remains homogeneous Neumann.
Also, $w(0)=w_0$ and $w(T)=0$.

Combining \eqref{bounded_interface_trace},
\eqref{internal_control_cost} and \eqref{reflected_data} gives
\eqref{cntrl_estmt}, after increasing the geometric constant.

If $w_0\in L^\infty$, then
$\|w^0\|_{L^\infty(Q)}\leq\|w_0\|_{L^\infty}$
by \eqref{homogeneous_maximum}.
In one space dimension $H^1(I)$ is continuously embedded in
$L^\infty(I)$. Hence \eqref{H2-estimate} gives
$$
\|w\|_{L^\infty(Q)}
\leq \|w_0\|_{L^\infty}
+C\|Y\|_{C([0,T];H^1(I))}
\leq \exp\!\big(C\kappa\big)\|w_0\|_{L^\infty}.
$$
Finally, \eqref{sln_estmt} and
$\|g\|_{L^2(0,T)}\leq T^{1/2}\|g\|_{L^\infty(0,T)}$
give the $V^1(Q)$ bound in \eqref{linear_controlled_state}.
\end{proof}

\section{Proof of the nonlinear result}\label{sc4}

We prove Theorem~\ref{main_thm} using the uniform linear estimates.
We fix the source in the elliptic equation, control the resulting
linear parabolic equation, and find a fixed point in $L^2(Q)$.

\subsection{Frozen systems and uniform bounds}

Fix $T>0$ and set
\begin{equation}\label{kappa_0}
\kappa_0=1+T+1/T,\qquad
\mathcal K=\{\xi\in L^2(Q):\|\xi\|_{L^\infty(Q)}\leq1\}.
\end{equation}
The set $\mathcal K$ is nonempty, closed and convex in $L^2(Q)$.
For $\xi\in\mathcal K$, let $u_\xi$ solve, for a.e.\ $t$,
\begin{equation}\label{fix_ellptc}
\begin{cases}
(u_\xi)_{xx}-u_\xi=\xi,&x\in(0,\mathrm L),\\
(u_\xi)_x(t,0)=(u_\xi)_x(t,\mathrm L)=0.
\end{cases}
\end{equation}
Proposition~\ref{prp21} gives
$u_\xi\in L^\infty(0,T;W^{2,p}(0,\mathrm L))$ for every
$1<p<\infty$, and $\|u_\xi\|_{L^\infty(Q)}\leq1$.
Define $B_\xi=f(u_\xi)$. Then
\begin{equation}\label{bnd_B}
\|B_\xi\|_{L^\infty(Q)}
\leq M_1:=\max_{|r|\leq1}|f(r)|.
\end{equation}
The frozen parabolic problem is
\begin{equation}\label{fix_prblic}
\begin{cases}
w_t=w_{xx}+B_\xi w_x,&(t,x)\in Q,\\
w_x(t,0)=g(t),\qquad w_x(t,\mathrm L)=0,&t\in(0,T),\\
w(0,x)=w_0(x),&x\in(0,\mathrm L).
\end{cases}
\end{equation}

By Theorem~\ref{thm31}, there is a constant $C_*>0$,
depending only on $\mathrm L$ and the bound $M_1$ in \eqref{bnd_B},
such that a null control
for \eqref{fix_prblic} can be chosen with
\begin{equation}\label{cntr_fixd}
\|g\|_{L^\infty(0,T)}
\leq A_T\|w_0\|_{L^\infty(0,\mathrm L)},
\qquad A_T:=\exp\!\big(C_*\kappa_0\big),
\end{equation}
and its state satisfies
\begin{equation}\label{sln_estmt_2}
\|w\|_{V^1(Q)}+\|w\|_{L^\infty(Q)}
\leq A_T\|w_0\|_{L^\infty(0,\mathrm L)}.
\end{equation}
We increase $C_*$ if necessary so that the same $A_T$ works in both
estimates. The constant $C_*$ is independent of $T$, $w_0$ and $\xi$.

\subsection{The set-valued map}

For each $\xi\in\mathcal K$, define
\begin{equation}\label{set}
\begin{split}
\Lambda(\xi):=\big\{w\in V^1(Q)\cap L^\infty(Q):\ &
\text{there is }g\in L^\infty(0,T)\text{ such that}\\
&\text{\eqref{fix_prblic}, \eqref{cntr_fixd},
and \eqref{sln_estmt_2} hold, and }w(T)=0\big\}.
\end{split}
\end{equation}
Thus
\begin{equation}\label{lmbda_map}
\Lambda:\mathcal K\longrightarrow 2^{L^2(Q)}.
\end{equation}
Its values are nonempty by the linear control theorem.
They are convex because the frozen equation is linear in $(w,g)$
and all the imposed bounds are convex.
Moreover, the full image $\Lambda(\mathcal K)$ is bounded in
$V^1(Q)$.

Lemma~\ref{lem:energy-compactness} now shows that
$\Lambda(\mathcal K)$ is relatively compact in $L^2(Q)$.
The bounds are uniform in $\xi$, so compactness holds for the whole
family of frozen systems. We include all controlled states satisfying
the bounds. This means that we do not need to choose a null control
continuously as the drift varies.

\subsection{Closed graph and terminal traces}

Let $\xi_n\to\xi$ strongly in $L^2(Q)$, with $\xi_n\in\mathcal K$,
and suppose that $w_n\in\Lambda(\xi_n)$ and $w_n\to w$ strongly
in $L^2(Q)$.
Choose corresponding controls $g_n$ satisfying \eqref{cntr_fixd}.
Passing to a subsequence, we have
\begin{equation}\label{cnvrgnt_1}
w_n\rightharpoonup w\ \text{in }V^1(Q),\qquad
g_n\rightharpoonup^* g\ \text{in }L^\infty(0,T).
\end{equation}
The estimates \eqref{cntr_fixd} and \eqref{sln_estmt_2}
pass to the limit by lower semicontinuity.

Corollary~\ref{cor:elliptic-time} gives
$$
\|u_{\xi_n}-u_\xi\|_{L^2(0,T;H^2(0,\mathrm L))}
\leq C\|\xi_n-\xi\|_{L^2(Q)}\longrightarrow0.
$$
Since $|u_{\xi_n}|,|u_\xi|\leq1$ and $f$ is Lipschitz on $[-1,1]$,
we have
\begin{equation}\label{B_n to B}
B_{\xi_n}\longrightarrow B_\xi
\quad\text{strongly in }L^2(0,T;L^\infty(0,\mathrm L)).
\end{equation}
After taking another subsequence, the convergence also holds pointwise
a.e.\ in $Q$. The uniform bound \eqref{bnd_B} then gives, for every
$\psi\in L^2(Q)$,
$$
\|(B_{\xi_n}-B_\xi)\psi\|_{L^2(Q)}\longrightarrow0
$$
by dominated convergence.
Combining this with the weak convergence of $(w_n)_x$ yields
\begin{equation}\label{weak_cnvgn_Bw}
B_{\xi_n}(w_n)_x\rightharpoonup B_\xi w_x
\quad\text{in }L^2(Q).
\end{equation}
The boundary traces of a test function in
$L^2(0,T;H^1(0,\mathrm L))$ belong to $L^2(0,T)$.
Since the time interval is finite, these traces are also integrable.
We can therefore pass to the limit in the boundary term using the
weak star convergence of $g_n$. Together with \eqref{weak_cnvgn_Bw},
this shows that $(w,g)$ solves \eqref{fix_prblic}.

For every fixed $t\in[0,T]$, evaluation at $t$ is a bounded linear
map from $V^1(Q)$ to $L^2(0,\mathrm L)$.
Consequently, \eqref{cnvrgnt_1} implies
$$
w_n(0)\rightharpoonup w(0),\qquad
w_n(T)\rightharpoonup w(T)
\quad\text{in }L^2(0,\mathrm L).
$$
Since $w_n(0)=w_0$ and $w_n(T)=0$, we obtain
$w(0)=w_0$ and $w(T)=0$.
Therefore $w\in\Lambda(\xi)$, and $\Lambda$ has a closed graph.
Together with the relative compactness of its full image, this
also proves that every value of $\Lambda$ is compact.

\subsection{Invariance and the fixed point}

\begin{proof}[Completion of the proof of Theorem~\ref{main_thm}]
Choose $c_1=C_*+1$ in \eqref{initil_cndtion}.
Then every $w\in\Lambda(\mathcal K)$ satisfies
$$
\|w\|_{L^\infty(Q)}
\leq A_T\|w_0\|_{L^\infty}
\leq \exp\!\big(-\kappa_0\big)\leq1.
$$
Hence $\Lambda(\mathcal K)\subset\mathcal K$.

To apply the fixed point theorem, we work on the closed convex hull
of the image of $\mathcal K$ in $L^2(Q)$. Define
$$
\mathcal D:=
\overline{\operatorname{co}}
\big(\overline{\Lambda(\mathcal K)}^{\,L^2(Q)}\big),
$$
where the outer closure is also taken in $L^2(Q)$.
The closure of $\Lambda(\mathcal K)$ is compact, and the closed
convex hull of a compact subset of a Banach space is compact.
A finite $\varepsilon$-net for the compact set also approximates every
convex combination by a combination of points in the net. These
combinations form a compact set, which proves total boundedness.
Thus $\mathcal D$ is nonempty, compact and convex.
Since $\mathcal K$ is closed and convex,
$\mathcal D\subset\mathcal K$.
For $\xi\in\mathcal D$ one has $\Lambda(\xi)\subset\mathcal D$.

The restriction of $\Lambda$ to $\mathcal D$ has nonempty compact
convex values and a closed graph. To see upper semicontinuity directly,
let $\mathcal O$ be an open set containing $\Lambda(\xi)$ for some
$\xi\in\mathcal D$. If upper semicontinuity failed, there would be
$\xi_n\to\xi$ in $\mathcal D$ and
$w_n\in\Lambda(\xi_n)\setminus\mathcal O$.
Compactness of $\mathcal D$ gives a subsequence with $w_n\to w$.
The closed graph implies $w\in\Lambda(\xi)\subset\mathcal O$,
contradicting $w_n\notin\mathcal O$ for all $n$.
The set-valued fixed point theorem
\cite{Glicksberg1952} gives $w\in\mathcal D$ such that
$w\in\Lambda(w)$.
Let $g$ be a corresponding control and put $u=u_w$.
Then \eqref{fix_ellptc} and \eqref{fix_prblic} give precisely
\eqref{equation_one}, and $w(T)=0$.
The regularity in the definition of weak solution follows from
the construction and Proposition~\ref{prp21}.

Finally, $w\in C([0,T];L^2(0,\mathrm L))$, and the elliptic solution
operator is bounded from $L^2(0,\mathrm L)$ to $H^2(0,\mathrm L)$.
Hence $u\in C([0,T];H^2(0,\mathrm L))$.
At $T$ its elliptic source is zero, so uniqueness in
\eqref{linr_ellptc} gives $u(T)=0$.
\end{proof}

\subsection{Uniqueness for a prescribed control}

We prove uniqueness of bounded weak solutions for a prescribed
initial datum and control.

\begin{prop}\label{prop:nonlinear-uniqueness}
For a prescribed $w_0$ and $g$, system \eqref{equation_one} has at most
one weak solution in the class of Definition~\ref{def:weak-solution}.
\end{prop}

\begin{proof}
Let $(w_1,u_1)$ and $(w_2,u_2)$ be two such solutions, and set
$W=w_1-w_2$ and $U=u_1-u_2$.
Choose $R$ so that $\|w_i\|_{L^\infty(Q)}\leq R$ for $i=1,2$.
The elliptic estimate gives $\|u_i\|_\infty\leq R$ and
$$
\|U(t)\|_\infty\leq C\|W(t)\|_2.
$$
By local Lipschitz continuity,
$\|f(u_1(t))-f(u_2(t))\|_\infty\leq C_R\|W(t)\|_2$.
Subtracting the parabolic equations gives
$$
W_t=W_{xx}+f(u_1)W_x+
\big(f(u_1)-f(u_2)\big)(w_2)_x,
$$
with zero initial datum and homogeneous Neumann data at both endpoints.
Testing with $W$ yields
$$
\frac12\frac{d}{dt}\|W\|_2^2+\|W_x\|_2^2
\leq C_R\|W_x\|_2\|W\|_2
+C_R\|(w_2)_x\|_2\|W\|_2^2.
$$
After Young's inequality,
$$
\frac{d}{dt}\|W\|_2^2+\|W_x\|_2^2
\leq C_R\big(1+\|(w_2)_x(t)\|_2\big)\|W\|_2^2.
$$
The coefficient on the right is integrable on $(0,T)$ because
$w_2\in V^1(Q)$. Gronwall's inequality and $W(0)=0$ give $W=0$.
Elliptic uniqueness then gives $U=0$.
\end{proof}

\begin{remark}
The proof uses only local Lipschitz continuity of $f$.
It does not assume $f(0)=0$, differentiability of $f$, positivity of
the states, or compatibility conditions on $w_0$ at the endpoints.
The one dimensional setting is used in the $L^6$ estimate for $Y_x$
and in the embedding $H^1(I)\hookrightarrow L^\infty(I)$.
\end{remark}

\begin{remark}
The Green representation for the elliptic equation in
\eqref{equation_one} shows that $w\geq0$ implies $u\leq0$.
The theorem therefore concerns signed states in this advection system.
No positivity assumption is imposed on either component.
\end{remark}

\section{Comments and open questions}\label{sc5}

We have proved local null controllability of \eqref{equation_one}
by a single bounded Neumann control. The elliptic equation also gives
$u(T)=0$. Lemma~\ref{lem:boundary-trace} allows a distributed control
on an extended interval to produce the boundary control.
We conclude with some consequences and possible extensions.

\subsection{Regularity and size of the initial data}

Theorem~\ref{thm31} already allows $L^2$ initial data for the linear
problem. In the nonlinear proof, we use the $L^\infty$ bound to keep
the controlled states in $\mathcal K$.
The elliptic solution operator maps $L^2(0,\mathrm L)$ into
$L^\infty(0,\mathrm L)$. It is therefore natural to ask whether the
nonlinear result can also be proved for small $L^2$ initial data.
This would require a different invariant set and a solution space
that allows unbounded initial data. In particular,
$w\in L^\infty(Q)$ cannot be required for every $w_0\in L^2$.

If $f$ is globally bounded, the same fixed point proof gives null
controllability for every $w_0\in L^\infty(0,\mathrm L)$.
We replace the unit source ball
by $\{\xi\in L^2(Q):\|\xi\|_{L^\infty(Q)}\leq R\}$.
The bound for $B_\xi$ is then independent of $R$. The linear estimates
bound the controlled state by $C_T\|w_0\|_\infty$, where $C_T$
depends only on $T$, $\mathrm L$ and $\|f\|_{L^\infty(\mathbb R)}$.
Choosing $R\geq C_T\|w_0\|_\infty$ keeps the controlled state in
the source ball, so the same fixed point argument applies.
For unbounded $f$, the drift bound increases with $R$, and this
argument may fail to give an invariant ball.

\subsection{Control cost and the smallness condition}

The bounds in \eqref{cntrl_estmt} and \eqref{initil_cndtion} are
sufficient, but we have not proved that their dependence on $T$ and
the drift bound is optimal. It would be useful to compare upper and
lower bounds for the boundary control, starting with constant drift.
Another question is how much of the smallness condition comes from
the linear control cost and how much comes from the nonlinear
invariant set.

\subsection{Higher dimensions and other boundary conditions}

An extension to higher dimensions would require a suitable vector
drift and a control on part of the boundary. Our proof uses
$H^1(I)\hookrightarrow L^\infty(I)$ and the $L^6$ estimate for $Y_x$.
Together, these estimates make the differentiated heat kernel integral
convergent. The corresponding energy estimates in higher dimensions
are not enough for the same bounded trace argument. Reflection across a curved
boundary also changes the equation, so new trace estimates are needed.
Robin conditions would require a different extension in place of the
even reflection used for homogeneous Neumann data.
The results in \cite{BhandariBoyer2021,BoyerOlive2024} treat related
boundary control problems. Further estimates would be needed to apply
their methods to \eqref{equation_one}.

\subsection{Constant states and nonconstant trajectories}

For every $c\in\mathbb R$, the pair $(w,u)=(c,-c)$ is a stationary
solution with $g=0$. Local exact controllability to this pair follows
from Theorem~\ref{main_thm}: the variables $w-c$ and $u+c$ satisfy
the same system with $f$ replaced by $f_c(s):=f(s-c)$.
The size of the admissible neighbourhood may depend on $c$.

Subtracting the equations for a nonconstant reference trajectory
$(\bar w,\bar u)$ gives the additional term
$[f(u)-f(\bar u)]\bar w_x$. The elliptic equation makes this term
nonlocal in the perturbation. Controlling the resulting equation
would require more than a bound for the drift.
The methods for controllability to trajectories and nonlocal couplings
in \cite{Guo_Zhu_Zhang,BhandariHernandez2025} may help in studying
this problem. The assumptions on the reference trajectory would
also have to be specified.

\subsection{Feedback stabilization}

The control constructed here depends on the prescribed initial datum
and control time. The fixed point argument does not give a feedback
law. One question is to construct a boundary feedback that
stabilizes the nonlinear system near zero. The linearization is the
heat equation with constant drift $f(0)$ and Neumann boundary
conditions. A feedback must control its constant mode, which does not
decay without control. If only partial measurements are available,
an estimate of the state is also needed to compute the boundary input.
The methods in \cite{Barbu18,coronbook07} may help with this
construction. For our system, the stability estimates must also
account for the nonlocal drift.

\section*{Acknowledgements}
The author acknowledges support from the Indian Institute of Technology
Bombay through an Institute Postdoctoral Fellowship (IPDF).

\bibliographystyle{elsarticle-num}
\bibliography{bib}
\end{document}